\documentstyle[11pt,leqno]{article}
\newtheorem{th}{Theorem}[section]
\newtheorem{lem}[th]{Lemma}
\newtheorem{pr}[th]{Proposition}
\newtheorem{cor}[th]{Corollary}

\title{\Large\bf{Conformally flat
Einstein-like 4-manifolds \\ and
conformally flat Riemannian 4-manifolds \\all of whose  Jacobi
operators have parallel \\ eigenspaces along every geodesic}}
\author{{\sc Stefan Ivanov} \thanks{The author supported
by Contract
MM 423/1994 with the Ministry of Science and Education of
Bulgaria and by
Contract 219/1994 with the University of Sofia "St. Kl.
Ohridski".} \hspace{5mm} {\sc Irina Petrova} \thanks{The author supported
by Contract
MM 413/1994 with the Ministry of Science and Education of
Bulgaria}}
\date{}
\begin{document}
\maketitle
\thispagestyle{empty}
\begin{center}
\end{center}
\begin{abstract}
A local classification of all locally conformal flat
Riemannian $4$-manifolds
whose Ricci tensor satisfies the equation $\nabla \left (
\rho-\frac{1}{6}sg\right )=\frac{1}{18}ds\odot g$ as well as a local
classification of all locally conformal flat Riemannian $4$-manifolds  for
which all Jacobi operators have parallel eigenspaces  along
every geodesic is
given. Non-trivial explicit examples are presented.  The
problem of local description of self-dual Einstein-like
$4$-manifolds is also treated. A complete explicit solution
of the St\"ackel system in dimension $4$ is obtained.
\\[20mm]
{\bf Running title:} Einstein-like 4-manifolds
\\[5mm]
\\[5mm]
${\bf 1991 \quad MS \quad Classification: } 53B20; 53C15; 53C55; 53B35$
\end{abstract}
\newpage
\section{Introduction}
One of the most important Riemannian metric on a smooth manifold is the
Einstein metric. In $1978$ A.Gray \cite{6} suggested to study three
generalizations of Einstein metrics. Background for his
investigations  was the
fact that any  Einstein metric has parallel Ricci tensor, and
conversely, that
any Riemannian manifold with parallel Ricci tensor is locally a product of
Einstein metrics. Using representation theory of the orthogonal group Gray
decomposed the covariant derivative of the Ricci tensor of a Riemannian
manifold into its irreducible components and derive so in a
natural way  three
classes of Riemannian manifolds, namely: the class $\cal B$ of
Riemannian manifolds whose Ricci tensor is a Codazzi tensor
(this is precisely the class of Riemannian manifolds with
harmonic curvature), the class $\cal U$
of Riemannian manifold with Killing Ricci tensor and the class $\cal Q$ of
Riemannian manifolds whose Ricci tensor satisfies
\begin{equation}\label{1}
   \nabla \left ( \rho-\frac{1}{2n-2}sg\right )=
   \frac{n-2}{2(n+2)(n-1)}ds\odot g
\end{equation}
where $n$ is the dimension of the manifold, $g$ its Riemannian metric,
$\nabla$ the Levi-Civita connection of the metric $g$, $\rho$ the
Ricci tensor, $s$ the scalar curvature and $\odot$ denotes the symmetric
product of symmetric tensors.
\par
In this paper we concentrate our attention on $\cal Q$-spaces
(for the sake of
brevity we call manifolds belonging to the class $\cal Q$ also $\cal
Q$-spaces).
\par
A complete local description (up to an isometry) of the
$3$-dimensional ${\cal Q}$-spaces  is given by J. Berndt in \cite{2}.
But, in higher dimensions
only few examples for $\cal Q$-spaces with non-parallel Ricci tensor are
known. A class of examples for $\cal Q$-spaces was presented in
\cite{Bes}, p.448. These examples arise as  certain bundles over
1-dimensional
manifolds whose fiber is an Einstein manifold with negative
scalar curvature.
If the fiber is a space of constant sectional curvature then
these examples are
locally conformal flat $\cal Q$-spaces. Besse considered also the
problem of local description of $\cal Q$-spaces, but he could not
obtain any classification of these spaces. As we know, this
problem is still open.
\par
In the paper we give a particular answer to this classification
problem.
We describe locally (up to an isometry) the 4-dimensional
locally conformal
flat $\cal Q$-spaces. Our considerations based on the fact
that for every locally conformal flat $\cal Q$-space  the
 Jacobi operator and its covariant derivative commute. Riemannian
manifolds for
which the Jacobi operator has this property are introduced by
J.Berndt and L.
Vanhecke in \cite{3} as a natural generalization of locally symmetric
Riemannian spaces.  These spaces are called ${\cal P}$-spaces in
\cite{3}. A
${\cal P}$-space is by definition a Riemannian manifold all of whose
Jacobi operators have parallel eigenspaces along every geodesic. The
${\cal P}$-spaces are also of special interest (see
\cite{4,5,BPV,btv}).
A complete local classification of the $2$-dimensional and the
$3$-dimensional
$\cal P$-spaces is given in \cite{3}. In higher dimensions however,  few
examples for not locally symmetric ${\cal P}$-spaces are known.
As we know, the
problem of local classification of $\cal P$-spaces is not yet
solved in higher dimensions.
\par
In this note we classify locally the $4$-dimensional locally conformal flat
$\cal Q$-spaces and locally conformal flat ${\cal P}$-spaces.
To this end, we write down $10$ kinds of explicit Riemannian metrics
on ${\bf R^4}$. Following the ideas of L.P.Eisenhart in
\cite{e2} we show that
these metrics give a complete solution of the so called St\"ackel system in
dimension $4$.
The present classification provides also some new examples of
$\cal Q$-spaces and some new examples of ${\cal P}$-spaces which seem to be
unknown up to now.  The purpose of this note is to prove the following two
local structure theorems
\begin{th}\label{1.1}
Let $({\bf M},g)$ be a connected $4$-dimensional locally
conformal flat $\cal P$-space.\\
Then $({\bf M},g)$ is locally (almost everywhere) isometric to one of the
following spaces
\par I) a real space form;
\par II) a Riemannian product of two $2$-dimensional real space forms of
opposite constant sectional curvatures;
\par $III_1$) a warped product $\bf B^1\times_f N^3$, where $\bf B^1$ is a
$1$-dimensional space, $\bf N^3$ is a $3$-dimensional real space
form and {\bf f} is a positive smooth function on ${\bf B^1}$.
\par IV) a warped product $\bf B^2\times_f N^2$, where $\bf N^2$ is a
$2$-dimensional
real space form of constant sectional curvature ${\cal K}_{\rm N}$, $\bf B^2$
is (a domain of) $\bf R^2$ with the following Riemannian metric
$$  g=\beta (x)\gamma (x)dx^2+\left( \left( \beta (x)y+
      \alpha (x)\right)dx+dy \right)^2,  $$
where
$$  \beta (x)=-\frac{{\cal K}'(x)}{2\left ({\cal K}(x)+
    \frac{{\cal K}_{\rm N}}{cA}+A\right )},\quad
    \gamma (x)=\frac{{\cal K}'(x)}{{\cal K}^2(x)-
               {\left (\frac{{\cal K}_{\rm N}}{cA}\right )}^2},  $$
$\alpha (x)$, ${\cal K}(x)$ are smooth functions of $x$, ${\cal K}'(x)\not=0$,
$c,A$ are constants,
different from zero, such that ${\cal K}(x)\not=
\pm\frac{{\cal K}_{\rm N}}{cA}$,\quad ${\cal K}(x)\not=
-\frac{{\cal K}_{\rm N}}{cA}-A$,\quad $c{\cal K}(x)+\frac{{\cal K}_{\rm N}}{A}
>0$\quad for every $x$. The positive function $f$ is determined by
$${\rm f}^2(x)=c{\cal K}(x)+\frac{{\cal K}_{\rm N}}{A}.$$
The Gauss curvature of $\bf B^2$ is equal to $\cal K$;
\par
V) a warped product $\bf B^2\times_f N^2$, where $\bf N^2$ is a
$2$-dimensional
real space form of constant sectional curvature ${\cal K}_{\rm N}$, $\bf B^2$
is (a domain of) $\bf R^2$ with the following Riemannian metric
$$  g=E(x)\left( \mu \left(x,y\right)+D\left(x\right)\right )dx^2+
      {\left (\frac{\mu_x(x,y)}{\mu_y(x,y)}dx+dy\right )}^2, $$
where
$$  E(x)=\frac{(D'(x))^2}{2D^2(x)\left ((D(x)-C)^2+\frac{2{\cal K}_{\rm N}}
         {cD(x)}+e\right )},
      $$
$D(x)$ is a smooth function of $x$, $D(x)\not=0$, $D'(x)\not=0$ for every
$x$, $C,c,e$ are constants, $c\not=0$ and $\mu(x,y)$ is a smooth solution of
the following PDE
$$  \mu_y^2(x,y)=\frac{2\mu(x,y)}{\mu(x,y)+D(x)}
    \left [\mu(x,y)\left ( (\mu(x,y)+C)^2+e\right )-\frac{2{\cal K}_{\rm N}}
    {c}\right ] $$
such that $\mu_y(x,y)\not=0$,\quad $cD(x)\mu (x,y)>0$,\quad $E(x)(\mu (x,y)+
D(x))>0$ for every $(x,y)$.
\\
The positive function $f$ is determined by
$${\rm f}^2(x,y)=cD(x)\mu (x,y).$$
The Gauss curvature ${\cal K}$ of $\bf B^2$ is equal to
${\cal K}=\mu-D+C$; \par
VI) a Riemannian manifold with the Riemannian metric of the form
$$
g=dx_1^2+\phi^2(x_1)dx_2^2+r\left(\phi^2\left(x_1\right)+a\right
)dx_3^2+q\left( \phi^2\left(x_1\right)+b\right)dx_4^2,
$$
where $a,b,q,r$ are positive constants, $a\not=b$ and
$\phi(x)$ is a non-constant solution of the following ordinary differential
equation
$$  (\phi'(x))^2=\phi^4(x)+(a+b)\phi^2(x)+ab. $$
\indent
VII) a Riemannian manifold with the Riemannian metric of the form
$$
g={{\Large \sigma }\atop{1,2,3,4}}\quad F_1(x_1)\vert
x_1-x_2\vert \vert x_1-x_3\vert \vert x_1-x_4\vert dx_1^2,
$$
where  ${\Large \sigma}$ denotes the cyclic sum and $1/F_i(x_i)=P(x_i)
(i=1,2,3,4)$ with a polynomial $P(x)$ of degree six.
\par
VIII) a warped product $\bf B^3\times_f N^1$, where $\bf N^1$ is a
one-dimensional Riemannian manifold, $\bf B^3$ is a three-dimensional
Riemannian manifold with Riemannian metric of the form
$$
g={{\Large \sigma }\atop{2,3,4}}\quad F_2(x_2)\vert x_2-x_3\vert
\vert x_2-x_4\vert dx_2^2, $$
where  ${\Large \sigma}$ denotes the cyclic sum and
$1/F_i(x_i)=P(x_i),
(i=2,3,4)$ with a polynomial $P(x)$ of degree five having  zero
as a root. The function $f$ is given by
$$f^2=\vert x_2x_3x_4 \vert .$$
\indent
IX) a Riemannian manifold with Riemannian metric of the form
$$
 g=(x_3-b)(x_4-b)dx_1^2+x_3x_4dx_2^2+{{\Large \sigma
}\atop{3,4}}\quad F_3(x_3)\vert x_3-x_4\vert dx_2^2, $$
where  ${\Large \sigma}$ denotes the cyclic sum and $1/F_i(x_i)=P(x_i)
(i=3,4)$ with a polynomial $P(x)$ of degree four given by
$P(x)=(x-b)(a_3x^3+a_2x^2+a_1x), \quad b\not=0$.
\par
Conversely, every Riemannian manifold of type I), II), $III_1$), IV), V),
VI), VII), VIII) and IX) is a locally conformal flat $\cal P$-space.
\end{th}
\begin{th}\label{1.2}
Let $({\bf M},g)$ be a connected $4$-dimensional locally
conformal flat $\cal
Q$-space. Then $({\bf M},g)$ is a $\cal P$-space and it is locally
(almost everywhere) isometric to one of the spaces I), II), IV),
V), VI), VII), VIII), IX) or to
\par
$III_2$) a warped product $\bf B^1\times_f N^3$ of a
$1$-dimensional base $\bf
B^1$ and a $3$-dimensional leaf $\bf N^3$ with constant sectional curvature
${\cal K}_{\rm N}$, and  $F=1/\rm f$ is a positive solution of the following
second order differential equation
\begin{equation}\label{w1}
  ({\rm F}''(x))^2=2{\cal K}_{\rm N}{\rm F}^3+c{\rm F},
\end{equation}
where $c$ is a constant.
\par
Conversely, every Riemannian manifold of type I), II), $III_2$), IV), V),
VI), VII), VIII) and IX) is a locally conformal flat $\cal Q$-space.
\end{th}
If the Jacobi operator has  pointwise constant eigenvalues then
the Riemannian manifold is said to be a pointwise Osserman manifold. By the
result of P.Gilkey, A.Swann and L.Vanhecke in \cite{GSV}, a 4-manifold is
pointwise Osserman manifold iff locally there exists a choice of orientation
such that the manifold is self-dual and Einstein.
It is well known (see e.g. \cite{Bes}, 16.4) that every $4$-dimensional $\cal
Q$-space has harmonic Weyl tensor. The results of J.P.Bourguignon
\cite{Bou9}, A.Derdzinski \cite{Der3}, A.Derdzinski and C.L.Chen
\cite{De-Ch} (see
also \cite{Bes}, 16.29) imply that any self-dual or anti-self-dual
$4$-dimensional $\cal Q$-space is either with zero Weyl
tensor or Einstein. Combining these results with Theorem 1.2
we obtain
\begin{cor}\label{co1}
Let ({\bf M},g) be a connected 4-dimensional self-dual (or anti-self-dual)
$\cal Q$-space. Then ({\bf M},g) is locally (almost everywhere) either
Einstein (and hence pointwise Osserman) or isometric to one of the
spaces described in Theorem 1.2.
\end{cor}
Thus, the description of self-dual
4-dimensional $\cal Q$-space is reduced to the difficult problem of local
classification of
pointwise Osserman 4-manifolds which is equivalent to the local description of
self-dual Einstein 4-manifolds (see the recent work of
N.Hitchin \cite{Hi} for the latter problem). \\[2mm]
{\bf Remarks}.
{\bf 1.} The spaces described in $III_2)$ and VII) of Theorem \ref{1.2} are
the locally conformal flat $\cal Q$-spaces presented in \cite{Bes}.
\par
{\bf 2.} The spaces described in $III_1)$ and $III_2)$ have parallel Ricci
tensor iff
the function $\rm f$ is either constant or $\rm f$ is determined
by one of the following three conditions:
\par
\par a) ${\rm f}=\varepsilon\sqrt {{\cal K}_{\rm N}}x+b$, where
$\varepsilon=\pm 1$, $b=const$, ${\cal K}_{\rm N}>0$,
\par b) ${\rm f}=Ce^{ax}+De^{-ax}$, where $a,C,D$ are constants,
${\cal K}_{\rm N}+4CDa^2=0$,
\par c) ${\rm f}=C\sin (ax)+D\cos (ax)$, where $a,C,D$ are constants and
${\cal K}_{\rm N}=a^2(C^2+D^2)$ if ${\cal K}_{\rm N}>0$.
\par
{\bf 3.} The spaces of type VI) are always with non-parallel Ricci tensor.
\par
{\bf 4.} The spaces of types VII), VIII) and IX) are also with non-parallel
Ricci
 tensor. If we consider the polynomial P(x) of degree less then
six in the case
 VII), P(x) of degree less then five in the case VIII), P(x) of
degree less then
 four in the case IX), then the Ricci tensor is parallel, but
this leads us to spaces of constant curvature.
\par
{\bf 5.} Comparing Theorem 1.1 and Theorem 1 of \cite{Der1} we
can conclude that every Riemannian 4-manifold with harmonic
curvature whose Ricci tensor
has at most two distinct eigenvalues is a
$\cal P$-space.
\par
{\bf 6.}  A locally conformal flat Riemannian manifold has
constant Ricci eigenvalues,  i.e. it is curvature homogeneous, iff it
is locally symmetric Riemannian manifold \cite{van},\cite{9}. Thus, we
obtain by Theorem
1.1 and Theorem 1.2 examples of locally conformal flat $\cal Q$-spaces
and examples of locally conformal flat $\cal
P$-spaces which are not even curvature homogeneous since locally conformal
flat curvature homogeneous Riemannian $4$-manifolds are exactly
the spaces with
constant Ricci  eigenvalues. All these examples are not also
pointwise Osserman manifolds by Corollary 1.3.

\section{Preliminaries}

Let $({\bf M},g)$ be an $n$-dimensional Riemannian manifold and $\nabla$
the Levi-Civita connection of the metric g. The curvature $R$ of $\nabla$
is defined by $R(X,Y)=[\nabla_X,\nabla_Y]-\nabla_{[X,Y]}$ for
every vector fields $X,Y$ on {\bf M}. We denote
by $T_p{\bf M}$ the tangential space at a point $p\in{\bf M}$.
\par
Let $X\in
T_p{\bf M}$. The Jacobi operator is defined by
$$  \lambda_X(Y)=R(Y,X)X, \qquad Y\in T_p{\bf M}.  $$
Let $\gamma(t)$ be a geodesic on ${\bf M}$ and $\dot{\gamma}$ denotes its
tangent vector field. We consider the family of smooth self-adjoint Jacobi
operators along $\gamma$ defined by $\lambda_{\dot{\gamma}}=R(X,
\dot{\gamma})\dot{\gamma}$ for every smooth vector field $X$ along $\gamma$.
A Riemannian manifold $({\bf M},g)$ is said to be a $\cal
P$-space if the
operators $\lambda_{\dot{\gamma}}$ have parallel eigenspaces
along every geodesic on ${\bf M}$ \cite{3}.
\par
A Riemannian
manifold $({\bf M},g)$ is said to be locally conformal flat
if around every point $p\in{\bf M}$ there exists a metric $\bar g$ which is
conformal to $g$ and $\bar g$ is flat.
By the Weyl theorem, an $n$-dimensional $(n\ge 4)$
Riemannian manifold is locally conformal flat iff the Weyl curvature
tensor $W$ vanishes (see e.g. \cite{Sch}) i.e. the curvature of
the metric $g$ has the following form
\begin{equation}\label{2}
   R(X,Y,Z,U)=
-\frac{s}{(n-1)(n-2)}\left(g\left(Y,Z)g(X,U)-g(X,Z)g(Y,U\right)\right)+
\end{equation}
$$+ \frac{1}{n-2}\left(\rho\left(Y,Z)g(X,U)-\rho(X,Z)g(Y,U)+
                  g(Y,Z)\rho(X,U)-g(X,Z)\rho(Y,U\right)\right),  $$
where $X,Y,Z,U \in T_p{\bf M}$, $p \in {\bf M}$. The following condition
also holds
\begin{equation}\label{3}
   (\nabla_X \rho)(Y,Z)-(\nabla_Y \rho)(X,Z)=
   \frac{1}{2(n-1)}\left(X\left(s)g(Y,Z)-Y(s)g(X,Z\right)\right).
\end{equation}

\section{Proof of Theorem \ref{1.1}}

The Ricci operator $Ric$ is defined by $g(Ric(X),Y)=\rho(X,Y),
X,Y\in T_p{\bf M}, p\in {\bf M}$.
In every point $p \in {\bf M}$ we consider the Ricci operator
$Ric$ as a linear
self-adjoint operator on the tangential space $T_p{\bf M}$. Let
$\Omega $ be
the subset of {\bf M} on which the number of distinct
eigenvalues of $Ric$ is
locally constant. This set is open and dense on {\bf M}. We can
choose $C^{\infty }$ eigenvalue functions  of
$Ric$ on $\Omega $, say $r_1,r_2,r_3,r_4$, in such a way such
that they form the
spectrum of $Ric$ at each point of $\Omega $ (see for example
\cite{10,11,3}). We fix a point $p \in \Omega $. Then there
exists a local orthonormal frame field $E_1,E_2,E_3,E_4$ on an
open connected neighborhood $V$ of $p$ such that
$$
Ric(E_i) = r_i E_i, \qquad i = 1,2,3,4.
$$
There exists an open connected neighborhood $U_p\subset V$ of
$p$ such that we have either $r_i=r_j$
or $r_i\not=r_j$ everywhere on $U_p$,  $i,j \in \{1,2,3,4\}$, $i\not=j$
\par
We fix the neighborhood $U_p$ and set
$$  (\nabla_i \rho)_{jk}=(\nabla_{E_i}\rho)(E_j,E_k),\quad
 \omega_{ij}^k=g(\nabla_{E_i} E_j, E_k), \quad i,j,k \in \{1,2,3,4\}.  $$
We have
\begin{equation}\label{4}
(\nabla_i \rho)_{jk}=\delta_{jk}E_i(r_i)+(r_j-r_k)\omega_{ij}^k,
\quad i,j,k \in \{1,2,3,4\},
\end{equation}
where $\delta_{jk}$ is the Kroneker's symbol.
\par
Further, unless otherwise stated, the Latin indices $i,j,k,l$
will stand for any distinct integers from the set $\{1,2,3,4\}$.
\par
We need the following technical result
\begin{pr}\label{2.1}
Four-dimensional locally conformal flat Riemannian manifold
$({\bf M},g)$ is a
$\cal P$-space iff for every point $p \in \Omega$  the following
conditions hold on $U_p$
\begin{equation}\label{5}
   (\nabla_i \rho)_{jk}=0,
\end{equation}
\begin{equation}\label{6}
   (\nabla_j \rho)_{jk}=(\nabla_i \rho)_{ik},
\end{equation}
\begin{equation}\label{7}
   E_i(r_k)=E_i(r_j),
\end{equation}
\begin{equation}\label{8}
   E_i(r_i)=3\left(E_i(r_k)-2(\nabla_k \rho)_{ki}\right),
\end{equation}
\begin{equation}\label{9}
   (r_i-r_l)\left (E_i(r_i)-E_i(r_k)-2(\nabla_k \rho)_{ki}\right )=0.
\end{equation}
\end{pr}
{\it Proof of Proposition \ref{2.1}}. Let $X \in T_p{\bf M}$, $p
\in {\bf M}$.
We consider the following skew symmetric operator $L_X$ on $T_p{\bf M}$
defined by $L_X=\lambda'_X\circ\lambda_X-\lambda_X\circ\lambda'_X$,
where $\lambda'_X$ is the covariant derivative of the Jacobi operator
defined by $\lambda'_X(Y)=(\nabla_XR)(Y,X,X), Y\in T_p{\bf M}$.
We know from \cite{3}  that $({\bf M},g)$ is a
$\cal P$-space iff the Jacobi operator and its covariant derivative
commute, i.e. $L_X=0$. Because of (\ref{2}) and (\ref{3}), the latter
equality is equivalent to
\begin{equation}\label{10}
   \left ( \rho(Y,Y)-\rho(Z,Z) \right)(\nabla_X \rho)(Y,Z)+
   \left ( (\nabla_X \rho)(Z,Z)-(\nabla_X \rho)(Y,Y)\right )\rho(Y,Z)+
\end{equation}
$$  +\rho(Y,U)(\nabla_X \rho)(Z,U)-\rho(Z,U)(\nabla_X \rho)(Y,U)=0, $$
where $X,Y,Z,U$ are orthonormal vectors of $T_p{\bf M}$, $p \in {\bf M}$.
\par
Let $({\bf M},g)$ be a $\cal P$-space. Taking $X=E_i$, $Y=E_j$, $Z=E_k$,
$U=E_l$,  we obtain from (\ref{10}) that
$$  (r_j-r_k)(\nabla_i \rho)_{jk}=0.  $$
The latter equality and (\ref{4}) imply (\ref{5}). Setting
$X=\cos\alpha E_i+\sin\alpha E_j$,\quad
$Y=\sin\alpha E_i-\cos\alpha E_j$,\quad $Z=\cos\beta
E_k+\sin\beta E_l$,\quad $U=\sin\beta E_k-
\cos\beta E_l$ into (\ref{10}), we obtain
\begin{equation}\label{11}
\sin2\alpha \left\{2\left(r_i\sin^2\alpha +r_j\cos^2\alpha
        -r_k\cos^2\beta -r_l\sin^2\beta \right) \right. .
\end{equation}
$$
\left. \left( \cos\beta \left((\nabla_i
\rho)_{ik}-(\nabla_j \rho)_{jk}\right) +
\sin\beta \left((\nabla_i \rho)_{il}-(\nabla_j
\rho)_{jl}\right) \right)\right. -
$$
$$
\left.  - (r_k-r_l)\sin2\beta \left(\sin\beta
\left(
(\nabla_i
\rho)_{ik}-(\nabla_j \rho)_{jk}\right)-\cos\beta
\left(
(\nabla_i \rho)_{il}-(\nabla_j
\rho)_{jl}\right)\right)\right\}=0.
$$
Replacing $-\beta $ by $\beta $ in (\ref{11})  and
adding the obtained equation to (\ref{11}), we obtain
$$  \sin2\alpha \cos\beta \left((r_i-r_j)\sin^2\alpha
+r_j-r_k\right)\left( \left(\nabla_i \rho)_{ik}-(\nabla_j
\rho\right)_{jk}\right)=0.
$$
The latter equality implies (\ref{6}).
\\
The formula (\ref{7}) follows from (\ref{6}) and
(\ref{3}). We derive
(\ref{8})  from (\ref{7}) and (\ref{3}).
\\
Let $A=(a^1,a^2,a^3)$, $B=(b^1,b^2,b^3)$, $C=(c^1,c^2,c^3)$ be orthonormal
vectors in $\bf R^3$. \\
We set
$$  X=\sum_{i=1}^3 a^iE_i,\quad Y=\sum_{i=1}^3 b^iE_i,\quad
    Z=\sum_{i=1}^3 c^iE_i,\quad U=E_4.  $$
Taking into account (\ref{5}), (\ref{6}), (\ref{7}) and
(\ref{8}), we
obtain  from (\ref{11}) by straightforward
computations that
\begin{equation}\label{12}
   \sum_{i=1}^3\left [ a^ib^ic^i\left (\sum_{j=1}^3
   \left( (b^j)^2-(c^j)^2\right )r_j \right)+
   a^i\left( (c^i)^2-(b^i)^2 \right) \left( \sum_{j=1}^3 b^jc^jr_j \right)
   \right ].
\end{equation}
$$
\left[E_i(r_i)-E_i(r_4)-2(\nabla_4 \rho)_{4i}\right]=0.
$$
We replace $\bar A=(-a^1,a^2,a^3)$, $\bar B=(b^1,-b^2,-b^3)$,
$\bar C=(c^1,-c^2,-c^3)$ by $A,B,C$ in (\ref{12}). Subtracting the
obtained equation from (\ref{12}), we get
$$  a^1\left [ \left (\sum_{j=1}^3 \left ( (b^j)^2-(c^j)^2 \right )r_j
    \right ) b^1c^1+
\left (\sum_{j=1}^3 b^jc^jr_j\right ) \left ( (c^1)^2-
 (b^1)^2\right) \right].
$$
$$
\left[E_1(r_1)-E_1(r_4)-2(\nabla_4 \rho )_{41}\right]=0.
$$
It is easy to see that the latter equality is equivalent to the equation
$$  (r_2-r_3)\left ( E_1(r_1)-E_1(r_4)-2(\nabla_4 \rho)_{41}\right )=0.  $$
Similarly we get all other equalities in (\ref{9}).
\\
Conversely, let the identities (\ref{5}), (\ref{6}), (\ref{7}), (\ref{8}) and
(\ref{9}) hold. The condition (\ref{10}) is equivalent to the equality
\begin{equation}\label{13}
   \sum_{s,p,r=1}^4 a^s \left [ b^pc^r\left ( \sum_{j=1}^4\left ( (b^j)^2-
   (c^j)^2\right )r_j\right )+(c^pc^r-b^pb^r)\left (\sum_{j=1}^4 b^jc^jr_j
   \right ) \right.
\end{equation}
$$  \left . +c^pd^r\left ( \sum_{j=1}^4 b^jd^jr_j \right )-b^pd^r\left (
    \sum_{j=1}^4 c^jd^jr_j\right )\right ](\nabla_s \rho)_{pr}=0  $$
for any orthonormal vectors $A=(a^1,a^2,a^3,a^4)$, $B=(b^1,b^2,b^3,b^4)$,
$C=(c^1,c^2,c^3,c^4)$, $D=(d^1,d^2,d^3,d^4)$ in $\bf R^4$.  We
check by direct
computations that the formulas (\ref{5}), (\ref{6}),
(\ref{7}), (\ref{8}) and (\ref{9}) imply (\ref{13}). \hfill {\bf Q.E.D.}
\par
Further, we have the following fife possibilities for the
eigenvalues of the Ricci tensor on $U_p$:
\par A)
$r_1=r_2=r_3=r_4$,
\par B)
$r_1=r_2=r_3\not=r_4$,
\par C)
$r_1=r_2\not=r_3=r_4$,
\par D)
$r_1=r_2\not=r_i$,\quad $i=3,4$,\quad $r_3\not=r_4$,
\par E)
$r_i\not=r_j$,\quad $i,j \in \{1,2,3,4\}$, $i\not=j$.
\\
{\bf The case A.} In this case $(U_p,g)$ is a real space form.
\\
{\bf The case B.} We get from Proposition \ref{2.1} that the following
equalities hold
\begin{equation}\label{14}
   \omega_{ij}^4=0,\quad{\rm for}\quad i,j \in \{1,2,3\}, i\not=j,
\end{equation}
\begin{equation}\label{15}
   \omega_{11}^4=\omega_{22}^4=\omega_{33}^4,
\end{equation}
\begin{equation}\label{16}
   \omega_{44}^1=\omega_{44}^2=\omega_{44}^3=0.
\end{equation}
Since $({\bf M},g)$ is locally conformal flat, we have
\begin{equation}\label{17}
   0=R(E_2,E_1,E_1,E_4)=E_2(\omega_{11}^4).
\end{equation}
We obtain similarly
\begin{equation}\label{18}
   E_3(\omega_{11}^4)=E_1(\omega_{11}^4)=0.
\end{equation}
We set ${\cal F}=span\{E_1,E_2,E_3\}$, ${\cal F}^{\perp}=span\{E_4\}$,
$\eta=\omega_{11}^4 E_4$. Let $\omega_{11}^4=0$ on $U_p$. Then $\cal F$ and
${\cal F}^{\perp}$ are autoparallel distributions. Hence, $(U_p,g)$ is
isometric to a Riemannian product $\bf B\times N$
of an $1$-dimensional space $B$ and a
$3$-dimensional Einstein space $N$. Hence, $N$ is
a real space form.
\par
Let there exist $q \in U_p$ such
that $\omega_{11}^4(q)\not=0$. It follows from (\ref{14}), (\ref{15}),
(\ref{16}),
(\ref{17}) and (\ref{18}) that we can apply a result of Hiepko (see
\cite{8}, p.211) to conclude that the space
$(U_p,g)$ is locally isometric,
around $q$, to a warped product of type $\bf B\times_f N$ with an
1-dimensional base $\bf B^1$ and a 3-dimensional leaf $N^3$
here ${\bf f}> 0$
is a smooth function on $\bf B^1$.  The Riemannian
metric $g^{{\bf B^1\times_f N^3 }}$ on $\bf B^1\times_f N^3$ is
given by
$g^{{\bf B^1\times_f N^3 }}=g^{{\bf B^1}}+f^2g^{{\bf N^3}},$
where $g^{{\bf B^1}}$ and $g^{{\bf N^3}}$ are the Riemannian
metrics on $\bf
B^1$ and $\bf N^3$, respectively. Using the formulae for the
Ricci eigenvalues
of the warped product $\bf B^1\times_f N^3$ (see \cite{1},
p.204-210), we have
\begin{equation}\label{19}
   r_i=\frac{\tilde \rho_i}{\tilde {\rm f}^2}-
   2{\left ( \frac{E_4(\tilde{\rm f})}{\tilde{\rm f}}\right )}^2-
   \frac{E_4^2(\tilde{\rm f})}{\tilde{\rm f}},\quad i=1,2,3,\qquad
   r_4=-3\frac{E_4^2(\tilde{\rm f})}{\tilde{\rm f}},
\end{equation}
where $\sim$ denotes the corresponding lift of a function on $\bf B^1$ or
$\bf N^3$ to $\bf B^1\times_f N^3$ and $\rho_i$,
$i=1,2,3$, are the
Ricci eigenvalues of $\bf N^3$. The conditions $r_1=r_2=r_3$ imply that
$\bf N^3$ is a $3$-dimensional Einstein space and hence it is a real space
form. Thus, $III_1)$ follows.
\\
{\bf The case C.} In this case (\ref{8}) and (\ref{9}) imply that all Ricci
eigenvalues are constant.   Then we have locally symmetric
Riemannian space since it is locally
conformal flat \cite{van},\cite{9} and  I) and II) follow.
\\
{\bf The case D.} The conditions (\ref{5}),
(\ref{6}), (\ref{7}), (\ref{8}) and (\ref{9})
of Proposition \ref{2.1} are equivalent to the
following equations
\begin{equation}\label{20}
   \omega_{ij}^k=0\quad {\rm for}\quad \{j,k\}\not=\{1,2\},
\end{equation}
\begin{equation}\label{21}
   \omega_{33}^k=\omega_{44}^k=0,\quad k=1,2,
\end{equation}
\begin{equation}\label{22}
   E_k(r_1)=E_k(r_3)=E_k(r_4)=0,\quad k=1,2,
\end{equation}
\begin{equation}\label{23}
   \omega_{11}^k=\omega_{22}^k,\quad k=3,4,
\end{equation}
\begin{equation}\label{24}
   E_k(r_k)=\frac{3}{2}E_k(r_1)=\frac{3}{2}E_k(r_j)=
\end{equation}
$$  =6(r_1-r_k)\omega_{11}^k=6(r_j-r_k)\omega_{jj}^k,\quad k,j \in \{3,4\},
    k\not=j.  $$
Since $(U_p,g)$ is locally conformal flat, we have
\begin{equation}\label{25}
   0=R(E_2,E_1,E_1,E_3)=E_2(\omega_{11}^3).
\end{equation}
We obtain analogously that
\begin{equation}\label{26}
   E_2(\omega_{11}^4)=0,\quad E_1(\omega_{11}^3)=E_1(\omega_{22}^3)=0,\quad
   E_1(\omega_{11}^4)=E_1(\omega_{22}^4)=0.
\end{equation}
We set
$$  {\cal F}=span\{E_1,E_2\},\quad {\cal F}^{\perp}=span\{ E_3,E_4\},\quad
    \eta=\omega_{11}^3E_3+\omega_{11}^4E_4. $$
It follows from (\ref{20}), (\ref{21}), (\ref{23}), (\ref{25}),
(\ref{26}) and
the  mentioned above result of Hiepko (\cite{8},p 211) that $(U_p,g)$ is
isometric to  a warped product of type $\bf B^2\times_f N^2$ with a
2-dimensional base $\bf B^2$ and a 2-dimensional leaf $N^2$
where ${\bf f}> 0$
is a smooth function on $\bf B^2$.  The Riemannian
metric $g^{{\bf B^2\times_f N^2 }}$ on $\bf B^2\times_f N^2$ is
given by
$g^{{\bf B^2\times_f N^2 }}=g^{{\bf B^2}}+f^2g^{{\bf N^2}},$
where $g^{{\bf B^2}}$ and $g^{{\bf N^2}}$ are the Riemannian
metrics on $\bf
B^2$ and $\bf N^2$, respectively.  Let $\sim$ denote a lift of a
vector field or a function on $\bf B^2$ or $\bf N^2$ to $\bf
B^2\times_f N^2$. Let $\cal
K$ and ${\cal K}_{\rm N}$ be the sectional curvatures of $\bf
B^2$ and $\bf
N^2$, respectively  and $X_1,X_2$ are orthonormal vector fields
on $\bf N^2$.
It follows by the well-known formulas of the warped product (see
e.g. \cite{1}, p.204-210) that the vector fields
$$  E_i=\frac{\tilde{X_i}}{\tilde{\rm f}},\quad i=1,2  $$
are orthonormal eigenvectors of the Ricci tensor
everywhere on $\bf B^2\times_f N^2$. The corresponding
Ricci eigenfunctions are given by
\begin{equation}\label{27}
   r_1=r_2=\frac{\widetilde{{\cal K}_{\rm N}}}{\tilde{\rm f}^2}-
   \frac{\widetilde{\Delta {\rm f}}}{\tilde{\rm f}}-
   \frac{g(\widetilde{grad({\rm f})},\widetilde{grad({\rm f})})}
   {\tilde{\rm f}^2},
\end{equation}
where $\Delta{\rm f}$ is the Laplacian of $\bf f$.
\\
There exist orthonormal smooth vector fields $X_3,X_4$ on $\bf
B^2$ such that the vector fields $E_i=\tilde{X_i},\quad i=3,4$
are also orthonormal eigenvectors of the Ricci
tensor everywhere
on $\bf B^2\times_f N^2$. The corresponding
Ricci eigenfunctions are determined by
\begin{equation}\label{28}
   r_3=\tilde{\cal K}-2\frac{H^{\rm
f}(\widetilde{X_3},\widetilde{X_3})}{\tilde{\rm
f}},\quad
   r_4=\tilde{\cal K}-2\frac{H^{\rm
f}(\widetilde{X_4},\widetilde{X_4})}{\tilde{\rm
f}}. \end {equation}
We also have
\begin{equation}\label{47}
H^{\rm f}(\widetilde{X_3},\widetilde{X_4})=0,
\end{equation}
where $H^{\rm
f}(\widetilde{X_3},\widetilde{X_4})=\widetilde{X_3}\widetilde{X_3}-
(\nabla_{\widetilde{X_3}}\widetilde{X_4}){\rm f}$  is the Hessian
of f. \begin{lem}\label{l1}
a) The warped product $\bf B^2\times_f N^2$ is locally
conformal flat iff the following equality holds
\begin{equation}\label{29}
    \frac{\widetilde{{\cal K}_{\rm N}}}{\tilde{\rm f}^2}+
    \frac{\widetilde{\Delta{\rm f}}}{\tilde{\rm f}}-
    \frac{g(\widetilde{grad({\rm f})},\widetilde{grad({\rm f})})}
    {\tilde{\rm f}^2}+\tilde{\cal K}=0.
\end{equation}
b) A locally conformal flat warped product
$\bf B^2\times_f N^2$ with $r_1=r_2\not=r_i$, $i=3,4$, $r_3\not=r_4$ is a
$\cal P$-space iff the condition (\ref{24}) holds.
\end{lem}
{\it Proof of the Lemma \ref{l1}}. The condition
a) follows  from the warped
product formulas (see e.g. \cite{1}, p.204-211)
and (\ref{2}) by direct computations. It is easy to verify that
   the equalities (\ref{20}), (\ref{21}) and (\ref{23}) always hold on
$\bf B^2\times_f N^2$. Formula
(\ref{29}) implies that if $\bf B^2\times_f N^2$ is a locally conformal flat
then \begin{equation}\label{30}
   \widetilde{{\cal K}_{\rm N}}={\cal K}_{\rm N}=const
\end{equation}
In this case (\ref{22}) also holds. We apply
Proposition \ref{2.1} to complete the proof of
Lemma \ref{l1}. \hfill {\bf Q.E.D.}
\par
We shall prove that the conditions (\ref{29}), (\ref{24}) of Lemma \ref{l1}
are equivalent to the following conditions
\begin{equation}\label{32}
    r_1=2\tilde{\cal K}+C,\quad r_3=2\tilde{\cal
K}-D(x_3)+C,\quad r_4=3\tilde{\cal K}+D(x_3),
\end{equation}
where $C$ is a constant and $D(x_3)$ is a smooth function of
$x_3$ such that $E_3(D)=-E_3(\tilde{\cal K})$,
\begin{equation}\label{33}
\frac{{\cal K}_{\rm N}}{\tilde{\rm f}^2}-{\left ( \frac{E_3(\tilde{\rm f})}
   {\tilde{\rm f}}\right )}^2-{\left ( \frac{E_4(\tilde{\rm f})}
   {\tilde{\rm f}}\right )}^2=\frac{\tilde{\cal K}+C}{2},
\end{equation}
\begin{equation}\label{34}
3E_3(\tilde{\cal K})=E_3(r_3)=\frac{3}{2}E_3(r_4)=6(r_1-r_3)\omega_{11}^3=
   6(r_4-r_3)\omega_{44}^3,
\end{equation}
\begin{equation}\label{35}
3E_4(\tilde{\cal K})=E_4(r_4)=\frac{3}{2}E_4(r_3)=6(r_1-r_4)\omega_{11}^4=
   6(r_3-r_4)\omega_{33}^4.
\end{equation}
We have from (\ref{27}), (\ref{28}) and (\ref{29}) that
\begin{equation}\label{36}
   r_3+r_4-r_1=3\tilde{\cal K}.
\end{equation}
We obtain from (\ref{24}) and (\ref{36}) the equalities
\begin{equation}\label{nov1}
\frac{2}{3}E_3(r_3)=E_3(r_4)=E_3(r_1)=2E_3(\tilde{\cal K}),
\end{equation}
\begin{equation}\label{nov2}
\frac{2}{3}E_4(r_4)=E_4(r_3)=E_4(r_1)=2E_4(\tilde{\cal K}).
\end{equation}
Then (\ref{nov1}) and (\ref{nov2}) imply (\ref{32}).
We get using the formula  $\Delta{\rm f}=H^{\rm
f}(X_3,X_3)+ H^{\rm f}(X_4,X_4)$ and (\ref{32}) that
$$
-\frac{\widetilde{\Delta{\rm f}}}{\tilde{\rm
f}}=\frac{3\tilde{\cal K}+C}{2}. $$
We obtain (\ref{33}) substituting the latter equality into (\ref{29}).
Conversely,
(\ref{32}), (\ref{33}), (\ref{34}) and (\ref{35}) imply (\ref{29}) and
(\ref{24}). The equalities (\ref{32}), (\ref{nov1}) and
(\ref{nov2}) imply that
$r_3=r_4$ iff $D$ and $\tilde{\cal K}$ are constants.
\par
It follows from the formulae (see \cite{1})
\begin{equation}\label{cd1}
\omega_{11}^3=-\frac{E_3(\tilde{\rm f})}{\tilde{\rm f}},\quad
    \omega_{11}^4=-\frac{E_4(\tilde{\rm f})}{\tilde{\rm f}},\quad
    \tilde{\rm f}\not=const
\end{equation}
that  there exists a
neighborhood $V_q$ of almost every point $q \in U_p$ such that one
of the following possibilities holds on $V_q$:
\par
D1) $\omega_{11}^3\not=0$, $\omega_{11}^4\not=0$ everywhere on $V_q$,
\par
D2) $\omega_{11}^3\not=0$, $\omega_{11}^4=0$ everywhere on $V_q$.
\\
We shall describe below the metric on $\bf B^2$ and the function
$f$ on $\bf B^2$ in the both cases.
\\
{\bf The case D1.} We obtain $E_4(\tilde{\cal K})\not=0$ and
$E_3(\tilde{\cal K})\not=0$ from (\ref{34}), (\ref{35}) and
(\ref{cd1}). We consider a local chart $(V_q,x)$ of $\bf B^2$
with local coordinate fields
$  \frac{\partial}{\partial x_3},\quad \frac{\partial}{\partial x_4}. $
We set
$$  a=g\left (\frac{\partial}{\partial x_3},
       \frac{\partial}{\partial x_4}\right ),\quad
    b^2=g\left ( \frac{\partial}{\partial x_3},
       \frac{\partial}{\partial x_3}\right ),\qquad
 ^{\rm B}\omega_{ii}^j=g(\nabla_{X_i} X_i,X_j),\quad i,j \in \{3,4\}.  $$
If $h$ is a smooth function on $\bf B^2$ we denote by $h,_i,h,_{ij}, \quad
i,j=3,4$ its first and second partial
derivatives with respect to the local coordinates $x_3,x_4$.
We have  (changing the sign of $X_3$ if it is necessary) that
$$
X_3=\frac{\frac{\partial}{\partial
x_3}-a\frac{\partial}{\partial x_4}}{\psi},
\quad X_4=\frac{\partial}{\partial x_4},
\quad {\rm where}\quad \psi=\sqrt{b^2-a^2}.
$$
Then, the metric $g$ on $\bf B^2$ is given by
\begin{equation}\label{ov1}
g=\psi^2dx_3^2+(adx_3+dx_4)^2.
\end{equation}
We calculate that
$$
^{\rm B}\omega_{44}^3=\frac{a,_4}{\psi},\quad
    ^{\rm B}\omega_{33}^4=-\frac{\psi,_4}{\psi}.
$$
\begin{equation}\label{nov3}
{\cal K}=X_3(^{\rm B}\omega_{44}^3)+X_4(^{\rm B}\omega_{33}^4)^2-
   (^{\rm B}\omega_{33}^4)^2-(^{\rm B}\omega_{44}^3)^2.
\end{equation}
Using the  notations in the case $D1$, we get
that the system of PDE (\ref{32}), (\ref{33}),
(\ref{34}) and (\ref{35}) is equivalent to the following system of PDE
\begin{equation}\label{37}
   \frac{{\rm f},_4}{\rm f}=\frac{{\cal K},_4}{2({\cal K}+D-C)},
\end{equation}
\begin{equation}\label{38}
   \frac{{\psi},_4}{\psi}=\frac{{\cal K},_4}{2({\cal K}+2D-C)},
\end{equation}
\begin{equation}\label{39}
   X_3({\cal K})=-X_3(D),\quad (i.e.\quad {\cal K},_3-a{\cal K},_4=-D,_3),
\end{equation}
\begin{equation}\label{40}
   \frac{X_3({\rm f})}{\rm f}=\frac{D,_3}{2D\psi},
\end{equation}
\begin{equation}\label{41}
   a,_4=-\frac{D,_3}{2({\cal K}+2D-C)},
\end{equation}
\begin{equation}\label{42}
   \frac{{\cal K}_{\rm N}}{{\rm f}^2}-
   {\left ( \frac{X_3({\rm f})}{\rm f}\right )}^2-
   {\left ( \frac{{\rm f},_4}{\rm f}\right )}^2=\frac{{\cal K}+C}{2},
\end{equation}
where $D(x_3)\not=0$,\quad $D,_3\not=0$, \quad ${\cal K}\not=C-2D$,\quad
${\cal K}\not=C-D$, \quad ${\cal K},_4\not=0$. \\
The equations (\ref{nov3}), (\ref{28}) and (\ref{27}) take the form
\begin{equation}\label{43}
   {\cal K}=X_3\left (\frac{a,_4}{\psi}\right )+
   \left (-\frac{{\psi},_4}{\psi}\right ),_4-
   {\left ( \frac{a,_4}{\psi}\right )}^2-
   {\left ( -\frac{{\psi},_4}{\psi}\right )}^2,
\end{equation}
\begin{equation}\label{44}
   {\cal K}-D+C=-2\left (\frac{X_3^2({\rm f})}{\rm f}+
                 \frac{\psi,_4}{\psi}.\frac{{\rm f},_4}{\rm f}\right ),
\end{equation}
\begin{equation}\label{45}
   2{\cal K}+D=-2\left ( \frac{{\rm f},_{44}}{\rm
f}-\frac{a,_4}{\psi}.    \frac{X_3({\rm f})}{\rm f}\right ),
\end{equation}
\begin{equation}\label{46}
2{\cal K}+C=\frac{{\cal K}_{\rm N}}{{\rm f}^2}-\frac{\Delta{\rm f}}{\rm f}-
               {\left(\frac{X_3({\rm f})}{\rm f}\right )}^2-
               {\left ( \frac{{\rm f},_4}{\rm f}\right )}^2.
\end{equation}
Moreover, the equation (\ref{47}) has to be also satisfied.
\par
We set $\mu={\cal K}+D-C$. We get from (\ref{37}) and (\ref{38}) that
\begin{equation}\label{nov4}
{\rm f}^2=H(x_3)\mu(x_3,x_4), \qquad \psi^2=E(x_3)(\mu(x_3,x_4)+D(x_3)),
\end{equation}
where $H(x_3)$ and $E(x_3)$ are smooth functions of $x_3$. Using
(\ref{nov4}), we obtain from (\ref{40}) that
\begin{equation}\label{nov5}
H(x_3)=cD(x_3), \quad c=const.
\end{equation}
We get using (\ref{nov4}), (\ref{nov5}) and (\ref{42}) that (\ref{44}) is
equivalent to
\begin{equation}\label{nov6}
E=\frac{D,_3^2}{2D^2\left ( (D-C)^2+\frac{2{\cal K}_{\rm N}}{cD}+e
      \right )}, \quad e=const.
\end{equation}
The equation (\ref{43}) follows from (\ref{39}), (\ref{nov5})
and (\ref{nov6}). The equation (\ref{39}) implies
\begin{equation}\label{ov2}
a=\frac{\mu,_3}{\mu,_4}.
\end{equation}
We obtain from (\ref{42}) and (\ref{nov6}) that
\begin{equation}\label{nov7}
 \mu,_4^2=-\frac{2\mu}{\mu+D}\left( \left ( \mu(\mu+C)^2+e\right )-
   \frac{2{\cal K}_{\rm N}}{c}\right )
\end{equation}
We replace (\ref{41}) by
$  \left (\frac{\mu,_3}{\mu,_4}\right ),_4=-\frac{D,_3}{2(\mu+D)}$.
The latter equality follows from (\ref{nov7}). It is easy to verify that
the equations (\ref{45}),(\ref{46}) and (\ref{47}) are also consequences of
(\ref{nov4})-(\ref{nov7}).
The equalities (\ref{nov4})-(\ref{nov7}) together with (\ref{ov1})
imply V) of Theorem 1.1. The spaces described in V) are
locally conformal flat $\cal P$-space with
$r_1=r_2\not=r_3\not=r_4\not=r_1$ by
Lemma \ref{l1}, (\ref{nov1}) and (\ref{nov2}). This completes the
considerations in the case {\bf D1}.
\\
{\bf The case D2.} We obtain $E_4(\tilde{\cal K})=0$ and
$E_3(\tilde{\cal
K})\not=0$ from (\ref{34}), (\ref{35}) and (\ref{cd1}). Then (\ref{32}),
(\ref{nov1}) and (\ref{nov2}) imply
$$
r_1=2{\cal K}+C, \quad r_3=3{\cal K}+A,\quad r_4=2{\cal K}-A+C,
$$
where $A$ is a constant. As in the previous case, we consider a
chart $(V_q,x)$
of $\bf B^2$ and keep the notations in {\bf D1}. We proceed
similarly as in the
case {\bf D1}. The system of PDE, corresponding to the system
(\ref{37})-(\ref{46}) in the case {\bf D1}, can be written in
the following way $$
   {\cal K},_4=0, \qquad {\rm f},_4=0, \qquad    \psi,_4=0,
$$
\begin{equation}\label{51}
   \frac{X_3({\rm f})}{\rm f}=\frac{{\cal K},_3}{2\psi({\cal K}+A-C)},
\end{equation}
\begin{equation}\label{52}
   a,_4=-\frac{{\cal K},_3}{2({\cal K}+2A-C)},
\end{equation}
\begin{equation}\label{53}
   \frac{{\cal K}_{\rm N}}{{\rm f}^2}-
   {\left (\frac{X_3({\rm f})}{\rm f}\right )}^2=
   \frac {{\cal K}+C}{2},
\end{equation}
\begin{equation}\label{54}
   {\cal K}=X_3\left( \frac{a,_4}{\psi}\right )-
   {\left (\frac{a,_4}{\psi}\right)}^2,
\end{equation}
\begin{equation}\label{55}
   {\cal K}-A+C=\frac{2a,_4}{\psi}.\frac{X_3({\rm f})}{\rm f},
\end{equation}
\begin{equation}\label{56}
   2{\cal K}+A=-2\frac{X_3^2({\rm f})}{\rm f},
\end{equation}
\begin{equation}\label{57}
   2{\cal K}+C=\frac{{\cal K}_{\rm N}}{{\rm f}^2}-
\frac{\Delta {\rm f}}{\rm f}-{\left ( \frac{X_3({\rm f})}{\rm f}\right)}^2,
\end{equation}
where ${\cal K},_3\not=0,\quad {\cal K}\not=C-2A,\quad {\cal
K}\not=C-A,\quad     A\not=0.$
\\
We get from (\ref{51}), (\ref{52}) and (\ref{55}) that
\begin{equation}\label{nov8}
\psi^2=-\frac{{\cal K},_3^2}{2({\cal K}+2A-C)({\cal K}^2-(A-C)^2)}.
\end{equation}
The equation (\ref{54}) is a consequence of (\ref{52}) and (\ref{nov8}). The
equation (\ref{51}) implies
\begin{equation}\label{nov9}
{\rm f}^2 =c(K+A-C), \quad c=const.\not=0.
\end{equation}
We obtain from (\ref{53}) and (\ref{nov9}) that
\begin{equation}\label{no1}
C=A-\frac{K_N}{cA}.
\end{equation}
Using (\ref{no1}), we obtain from (\ref{52}) that
\begin{equation}\label{no2}
  a=-\frac{{\cal K},_3}{2\left ( {\cal K}+\frac{{\cal K}_{\rm N}}
      {cA}+A\right )}x_4+\alpha (x_3),
\end{equation}
where $\alpha(x_3)$ is a smooth function of $x_3$.
The equations (\ref{56}),(\ref{57}) and (\ref{47}) are also consequences of
(\ref{nov8})-(\ref{no2}).
\par
We substitute (\ref{no1}) into (\ref{nov8}) and
(\ref{nov9}).  Substituting
the obtained equations and (\ref{no2}) into (\ref{ov1}), we
get the spaces described in IV) of Theorem 1.1.
These spaces are locally conformal flat $\cal P$-space with
$r_1=r_2\not=r_3\not=r_4\not=r_1$
by Lemma \ref{l1}, (\ref{nov1}) and (\ref{nov2}). This completes
the considerations in the case {\bf D2}.
\\
{\bf The case E.} In this case, the formulae (\ref{5}), (\ref{6}), (\ref{7}),
(\ref{8}) and (\ref{9}) of Proposition \ref{2.1} are equivalent to
\begin{equation}\label{59}
   \omega_{ij}^k=0,
\end{equation}
\begin{equation}\label{60}
   (r_j-r_k)\omega_{jj}^k=(r_i-r_k)\omega_{ii}^k,
\end{equation}
\begin{equation}\label{61}
   E_k(r_i)=E_k(r_j),
\end{equation}
\begin{equation}\label{62}
   E_k(r_k)=\frac{3}{2}E_k(r_i),
\end{equation}
\begin{equation}\label{63}
   E_k(r_i)=4(r_i-r_k)\omega_{ii}^k.
\end{equation}
Let ${\cal F}_i$ be the orthonormal complement to $E_i$,
$i=1,2,3,4$. The equality (\ref{59}) implies that the distributions
${\cal F}_i$ are integrable. Then (decreasing, if it is necessary, the
neighborhood
$U_p$ of the point $p$) we can find a chart $(U_p,x)$ of $({\bf M},g)$,
(see \cite{7,e2}) such
that the metric $g$ is given by
\begin{equation}\label{64}
   g=\mu_1^2dx_1^2+\mu_2^2dx_2^2+\mu_3^2dx_3^2+\mu_4^2dx_4^2,
\end{equation}
where the functions $\mu_i=\mu_i(x_1,x_2,x_3,x_4)$ are smooth and strictly
positive and we have
\begin{equation}\label{65}
   X_i:=\frac{\partial}{\partial x_i}=\mu_i E_i,\quad i=1,2,3,4.
\end{equation}
We set
$$  \nu_s=\ln (\mu_s),\quad \nu_{s,t}=\frac{\partial}{\partial x_t}\nu_s=
    X_t(\nu_s),\quad \nu_{s,tr}=\frac{\partial^2}{\partial x_t \partial x_r}
    \nu_s, $$
$$  \nabla_{X_s}X_t=\sum_{r=1}^4 \Gamma_{st}^r X_r,\qquad
    R(X_t,X_r,X_s)=\sum_{u=1}^4R_{str}^u X_u  $$
We calculate using (\ref{64}) that
\begin{equation}\label{66}
   \Gamma_{ij}^k=0,\quad \Gamma_{ik}^k=\Gamma_{ki}^k=\nu_{k,i},\quad
   \Gamma_{ii}^k=-\frac{\mu_i^2}{\mu_k^2}\nu_{i,k},\quad
   \Gamma_{ii}^i=\nu_{i,i},
\end{equation}
\begin{equation}\label{67}
R_{jij}^i=\mu_j^2\left [\frac
{1}{\mu_j^2}(\nu_{i,j}\nu_{j,j}-\nu_{i,j}^2-
\nu_{i,jj})+\frac{1}{\mu_i^2}(\nu_{j,i}\nu_{i,i}-\nu_{j,i}^2-\nu
_{j,ii})\right .
\end{equation}
$$  \left . -\frac{1}{\mu_k^2}\nu_{i,k}\nu_{j,k}-
    \frac{1}{\mu_l^2}\nu_{i,l}\nu_{j,l}\right ],   $$
\begin{equation}\label{68}
  R_{kij}^i=-\nu_{i,jk}+\nu_{i,j}\nu_{j,k}+\nu_{i,k}\nu_{k,j}-
    \nu_{i,j}\nu_{i,k}.
\end{equation}
\begin{lem}\label{l2}
If the formulas (\ref{59})-(\ref{63}) are valid then the functions
$\nu_1,\nu_2,\nu_3,\nu_4$ satisfy the following system of non-linear PDE on
$U_p$.
\begin{equation}\label{69}
   \nu_{i,j}\nu_{j,k}+\nu_{i,k}\nu_{k,j}-\nu_{i,j}\nu_{i,k}=0,
\end{equation}
\begin{equation}\label{70}
   \nu_{i,jk}=0,
\end{equation}
\begin{equation}\label{71}
   \nu_{i,ij}+2\nu_{i,j}\nu_{j,i}=0.
\end{equation}
\end{lem}
{\it Proof of the Lemma}:
It follows from (\ref{65}) and (\ref{66}) that
\begin{equation}\label{72}
   \nu_{i,k}+\mu_k\omega_{ii}^k=0
\end{equation}
Then, the equality (\ref{60}) is equivalent to
\begin{equation}\label{73}
   (r_i-r_k)\nu_{i,k}=(r_j-r_k)\nu_{j,k}.
\end{equation}
If $\nu_{i,k}=0$ on $U_p$, we have from (\ref{72}) that $\nu_{j,k}=0$ which
implies (\ref{69}).
\\
Let there exist a point $q \in U_p$ such that
$\nu_{i,k}(q)\not=0$. The formula (\ref{69}) follows from (\ref{73})
and the equality
$$  \frac{r_i-r_k}{r_j-r_k}-1=-\frac{r_i-r_j}{r_k-r_j}. $$
\\
We obtain $R_{kij}^i=0$ from (\ref{2}) and (\ref{65}). The
latter equality, together with (\ref{68}) and (\ref{69}), implies (\ref{70}).
\\
If $\nu_{i,j}=0$ on $U_p$ then (\ref{71}) is valid. Let $\nu_{i,j}\not=0$. It
follows from (\ref{73}) that $\nu_{k,j}\not=0$ and the following equalities
hold
\begin{equation}\label{74}
   \frac{r_k-r_j}{r_i-r_j}=\frac{\nu_{i,j}}{\nu_{k,j}}.
\end{equation}
We have from (\ref{62}), (\ref{63}), (\ref{65}) and (\ref{72}) that
\begin{equation}\label{75}
   X_i(r_j-r_i)=-\frac{1}{2}X_i(r_j)+2(r_j-r_i)\nu_{j,i}.
\end{equation}
Using (\ref{61}), (\ref{74}) and (\ref{75}), we obtain consequently
$$  X_i(r_j-r_i)=X_i(r_k-r_i)=(r_j-r_i)X_i\left (
\frac{r_k-r_i}{r_j-r_i}\right )+
    \frac{r_k-r_i}{r_j-r_i}X_i(r_j-r_i),  $$
$$  \left ( 1-\frac{r_k-r_i}{r_j-r_i}\right )X_i(r_j-r_i)=
    (r_j-r_i)X_i\left (\frac {r_k-r_i}{r_j-r_i}\right ), $$
$$  2\left (1-\frac{r_k-r_i}{r_j-r_i}\right )(r_j-r_i)\nu_{j,i}=
    -(r_j-r_i)X_i\left ( 1- \frac{r_k-r_i}{r_j-r_i}\right ),  $$
$$2\frac{r_j-r_k}{r_j-r_i}=-X_i\left ( \frac{r_j-r_k}{r_j-r_i}\right ), $$
$$  2\frac{\nu_{i,j}}{\nu_{k,j}}\nu_{j,i}=-X_i\left ( \frac{\nu_{i,j}}
    {\nu_{k,j}}\right )  $$
The latter equality together with (\ref{70}) implies (\ref{71}).
So the whole Lemma is proved. \hfill {\bf Q.E.D.}
\par
The system (\ref{69})-(\ref{71}) of non-linear PDE is well known
as a system of
St\"ackel: it arises from a quantum mechanical problems; on the
other hand it
is shown in \cite{3} that this system describes 3-dimensional
${\cal P}$-spaces
with distinct Ricci eigenvalues; it is shown in \cite{IP1} that the local
description of Riemannian 3-manifolds with constant eigenvalues
of the natural
skew-symmetric curvature operator on each unit circle also reduces to this
system. The St\"ackel system is solved in Euclidean $3$-space by Weinacht
in \cite{Wei}, W.Blaschke in \cite{WB}.
L.P.Eisenhart solved in \cite{e2} this system completely in dimension
three and obtained solutions in higher dimensions. In \cite{e1} L.P.Eisenhart
classified all three dimensional St\"ackel systems in locally conformal flat
3-space. In order to complete the proof of Theorem 1.1 we have to solve  the
St\"ackel system in dimension $4$. Following the ideas of \cite{e2} we obtain
a complete solution of the St\"ackel system in dimension $4$, namely we have
\begin{lem}\label{l3}
The St\"ackel system (\ref{69})-(\ref{71}) has the following $10$ kinds of
solutions in dimension $4$:
\par
${\bf I_1.} \quad \mu _1^2=1, \mu _2^2=\eta (x_1), \mu _3^2=\eta (x_1)\psi
(x_2),
\mu_4^2=\eta (x_1)\phi (x_2)$, where $\eta '\not=0$, $\eta $ is a function of
$x_1$, $\psi $ and $\phi $ are   functions of $x_2$.
\par
${\bf I_2.} \quad \mu _1^2=1, \mu _2^2=\phi (x_1)(\xi (x_2)-\zeta
(x_3))(\xi
(x_2)-\eta (x_4)), \mu _3^2=\phi (x_1)(\zeta (x_3)-\xi
(x_2))(\zeta (x_3)-\eta (x_4)), \mu _4^2=\phi (x_1)(\eta
(x_4)-\xi (x_2))(\eta (x_4)-\zeta (x_3))$, where $\psi ,\phi ,
\eta,\zeta$ are functions of one variable and $\zeta '\not=0,\eta '\not=0$.
\par
${\bf I_3.} \quad \mu _1^2=\eta (x_1), \mu _2^2=x_1\xi (x_2), \mu _3^2=\mu
_4^2=x_1x_2(\phi (x_3)+\psi (x_4))$,
where $\eta ,\xi ,\psi , \phi $ are  functions of one variable.
\par
${\bf I_4.} \quad \mu _1^2=\mu _2^2=\phi (x_1)+\psi (x_2), \mu _3^2=\mu
_4^2=\xi (x_3)+\eta (x_4)$, where $\phi,\psi,\xi,\eta $ are
functions of one variable.
\par
${\bf I_5.} \quad \mu_1^2=\mu_2^2=(\xi(x_1)-\eta(x_2)), \mu_3^2=\mu_4^2
=\xi(x_1)\eta(x_2)(\psi(x_3)+\phi(x_4))$, where $\phi ,\psi ,\xi ,\eta $ are
functions of one variable and $\xi'\not=0, \eta'\not=0$.
\par
${\bf I_6.} \quad \mu_1^2=\phi(x_1), \mu_2^2=\psi(x_1),
\mu_3^2=\mu_4^2=x_1(\xi(x_3)+\eta(x_4))$, where $\phi ,\psi ,\xi ,\eta $ are
functions of one variable and $\phi >0, \psi >0$
\par
${\bf I_7.} \quad \mu _1^2=1, \mu _2^2=\phi (x_1), \mu _3^2=\psi (x_1), \mu
_4^2=\eta (x_1)$, where $\eta, \psi, \phi $ are  functions of $x_1$.
\par
${\bf I_8.} \quad \mu _i^2=\phi
_i(x_i)\vert x_i-x_j\vert \vert x_i-x_k\vert \vert x_i-v_l\vert $
for $i,j,k,l \in \{1,2,3,4\}$ distinct, where $\phi _i$ is a
function of $x_i$. \par
${\bf I_9.} \quad \mu _1^2=x_2x_3x_4, \mu _i^2=\phi_i (x_i)\vert
x_i-x_j\vert \vert
x_i-x_k\vert $ for $i,j,k \in \{2,3,4\}$ distinct, where
$\phi_i$ is a function of $x_i$.
\par
${\bf I_{10}.} \quad \mu _1^2=\vert x_3-b\vert \vert x_4-b\vert ,
\mu_2^2=x_3x_4,
\mu_i^2=\phi_i(x_i)\vert x_i-x_j\vert $ for $i,j \in \{3,4\}$
distinct, where $\phi_i$ is a function of $x_i$.
\end{lem}
{\it Proof of the Lemma}. We follow the considerations of
L.P.Eisenhart. In
\cite{e2} he showed how one could obtain all solutions of the
St\"ackel system
in an arbitrary dimension and he found four kinds of explicit solutions in
dimension three. We apply his considerations in \cite{e2},p.$291-294$ to the
four dimensional case. The formula $(3.14)$ of \cite{e2} implies the case
${\bf I_8}$. All solutions are given by formula $(3.9)$ of
\cite{e2}.  We write
down explicitly the conditions $(3.10)$ and $(3.11)$ of
\cite{e2} for four distinct
indices. All possible cases are given by the formulas $(3.17)-(3.19)$ of
\cite{e2}. We examine carefully each of these cases and we get
by a lengthy but
straitforward computations that all the remaining solutions of the St\"ackel
system in dimension four can be described by ${\bf I_1}-{\bf I_7},
{\bf I_9},{\bf I_{10}}$,
using a suitable changes of the local coordinates. \hfill {\bf Q.E.D.}
\par
To complete the proof of Theorem 1.1 we have to pick out those of the spaces
${\bf I_1}-{\bf I_{10}}$ which are locally conformal flat $\cal P$-spaces and
have four distinct eigenvalues of the Ricci operator. We have
\begin{lem}\label{nl4}
Let $(\bf M,g)$ be a $4$-dimensional Riemannian manifold with the
metric of the
form (\ref{64}) for which the conditions of Lemma \ref{l2} hold. Then
\par
i) $({\bf M},g)$ is locally conformal flat iff
\begin{equation}\label{d1}
R(E_i,E_l,E_l,E_i) +
R(E_k,E_j,E_j,E_k)=R(E_i,E_k,E_k,E_i)+R(E_j,E_l,E_l,E_j).
\end{equation}
\indent ii) If ({\bf M},g) is locally conformal flat with four
distinct Ricci eigenvalues then it is a $\cal P$-space iff
\begin{equation}\label{p1}
\frac{2}{3}X_i(r_i)=X_j(r_j)=4(r_i-r_j)\nu_{j,i}.
\end{equation}
\end{lem}
{\it Proof of the Lemma}. The proof of i) follows from
(\ref{2}) and (\ref{64})
by direct computations. The proof of ii) is a consequence of
(\ref{59})-(\ref{63}) and the conditions of the Lemma. \hfill
{\bf Q.E.D.}
\par
{\bf Spaces of type} ${\bf I_1,I_2}$  and ${\bf I_3}$. A space
of each of these three types is a warped product of a
$1$-dimensional base $B$ and a $3$-dimensional leaf $N$. Since
it is locally
conformal flat then the leaf $N$ is a space of constant
sectional curvature.
Then (\ref{19}) implies that the number of distinct Ricci
eigenvalues is at most two.
\par
{\bf Spaces of type} ${\bf I_4,I_5,I_6}$ and ${\bf I_{10}}$ {\bf
with $b=0$}. A space of any of these four types can be
regarded as a warped (or Riemannian) product of two Riemannian
surfaces. It follows by the
formulas for the Ricci tensor of a warped product (\ref{27})
that the number of distinct Ricci eigenvalues is at most three.
\par
{\bf Spaces of type} ${\bf I_7}$. Using (\ref{65}) and
(\ref{67}), we calculate consequently
\begin{equation}\label{83}
   R(E_1,E_j,E_j,E_1)=-\nu_{j,1}^2-\nu_{j,11}, \quad
R(E_i,E_j,E_j,E_i)=-\nu_{i,1}\nu_{j,1}
\end{equation}
for $i,j\in\{2,3,4\}$.
\begin{equation}\label{84}
r_1=-\nu_{2,1}^2-\nu_{3,1}^2-\nu_{4,1}^2-\nu_{2,11}-\nu_{3,11}-\nu_{4,11},
\end{equation}
\begin{equation}\label{85}
   r_i=-\nu_{i,11}-\nu_{i,1}(\nu_{2,1}+\nu_{3,1}+\nu_{4,1}),\quad i=2,3,4,
\end{equation}
Because of (\ref{d1}), we get from (\ref{85}) that
\begin{equation}\label{86}
   r_i-r_j=2\nu_{k,1}(\nu_{j,1}-\nu_{i,1}),\quad i,j,k \in \{2,3,4\}.
\end{equation}
It follows from (\ref{61}) that $X_1(r_i-r_j)=0$. The latter equality and
(\ref{86}) imply
\begin{equation}\label{87}
   (\nu_{j,11}-\nu_{i,11})\nu_{k,1}+(\nu_{j,1}-\nu_{i,1})\nu_{k,11}=0.
\end{equation}
We obtain using (\ref{83}) that the equality (\ref{d1}) for $l=1$
is equivalent to
\begin{equation}\label{88}
\nu_{j,11}-\nu_{i,11}=(\nu_{j,1}-\nu_{i,1})(\nu_{k,1}-\nu_{j,1}-\nu_{i,1})
\end{equation}
We get from (\ref{86}), (\ref{87}) and (\ref{88}) that
\begin{equation}\label{89}
   \nu_{i,1}\not=0,\quad \nu_{i,1}\not=\nu_{j,1},
\end{equation}
\begin{equation}\label{90}
   (\nu_{i,1}+\nu_{j,1}-\nu_{k,1})\nu_{k,1}-\nu_{k,11}=0,
\end{equation}
where $i,j,k$ belong to the set $\{2,3,4\}$.
\par
To complete the considerations in the case ${\bf I_7}$ we have to solve the
system of PDE (\ref{89}) and (\ref{90}). By direct computations
it follows that the conditions (\ref{90}) are equivalent to the
following equations
\begin{equation}\label{91}
   (\nu_{2,1}\nu_{3,1}),_1=(\nu_{2,1}\nu_{4,1}),_1=(\nu_{3,1}\nu_{4,1}),_1=
   2\nu_{2,1}\nu_{3,1}\nu_{4,1}.
\end{equation}
The latter equalities are equivalent to
\begin{equation}\label{92}
   \nu_{3,1}(\nu_{2,1}-\nu_{4,1})=D,\quad \nu_{4,1}(\nu_{4,1}-\nu_{3,1})=F,
   \quad \nu_{3,1}\nu_{4,1}=Pe^{2\nu_2}
\end{equation}
where $D$,$F$, $P$ are constants. We obtain from the symmetries
of (\ref{90}) that
\begin{equation}\label{93}
   \nu_{2,1}\nu_{3,1}=Qe^{2\nu_4},\quad \nu_{2,1}\nu_{4,1}=Re^{2\nu_3},
\end{equation}
for some constants $Q,R$. Then we have
$$  Qe^{2\nu_4}-Pe^{2\nu_2}=D,\quad Re^{2\nu_3}-Pe^{2\nu_2}=F,  $$
\begin {equation}\label{94}
   \nu_3=\frac{1}{2}\ln \left ( \frac{Pe^{2\nu_2}+F}{R}\right ),\quad
   \nu_4=\frac{1}{2}\ln \left ( \frac{Pe^{2\nu_2}+D}{Q}\right ),
\end{equation}
$$  \nu_{3,1}=\frac{Pe^{2\nu_2}\nu_{2,1}}{Pe^{2\nu_2}+F},\quad
    \nu_{4,1}=\frac{Pe^{2\nu_2}\nu_{2,1}}{Pe^{2\nu_2}+D}.  $$
Substituting the latter two equalities in the third equation of
(\ref{92}), we obtain
\begin{equation}\label{95}
   \nu_{2,1}^2=\frac{(Pe^{2\nu_2}+F)(Pe^{2\nu_2}+D)}{Pe^{2\nu_2}}.
\end{equation}
It follows from (\ref{89}), (\ref{92}) and (\ref{93}) that
\begin{equation}\label{96}
   F\not=D,\quad FDPQR\not=0.
\end{equation}
Conversely,  the formulas  (\ref{94}),
(\ref{95}) and (\ref{96}) imply (\ref{92}), (\ref{93}) and (\ref{89}).
\\
If we set
$  a=\frac{F}{P},\quad b=\frac{D}{P},\quad r=\frac{P}{R},\quad
    q=\frac{P}{Q}, $
than we obtain that (\ref{94}), (\ref{95}) and (\ref{96}) are equivalent to
$  \mu_3^2=r(\mu_2^2+a),\quad \mu_4^2=q(\mu_2^2+b),\quad
    \mu_{2,1}^2=(\mu_2^2+a)(\mu_2^2+b),  \quad
r\not=0,\quad q\not=0,\quad a\not=0,\quad b\not=0,\quad
a\not=b$. Setting $\phi =\mu_2$ we get the metrics described in the
case VI). \par
Conversely, for any metrics described in VI), we get from
(\ref{83}) that (\ref{d1}) is satisfied. We calculate for the Ricci
eigenvalues using (\ref{84}) and (\ref{85}) that
\begin{equation}\label{nov}
r_1=-6\phi^2-2a-2b, \quad r_2=-4\phi^2-2a-2b,
\end{equation}
$$r_3=-4\phi^2-2b, \quad r_4=-4\phi^2-2a.
$$
The equalities (\ref{nov}) imply that all the Ricci eigenvalues are distinct
and (\ref{p1}) is satisfied. Then Lemma \ref{nl4} proves VI).
\par
{\bf Spaces of type} ${\bf I_8}$. We set $\psi_s=\frac{1}{\phi_s(x_s)},
s=1,2,3,4$. We
multiply (\ref{d1}) by $\frac{1}{(x_i-x_l)(x_j-x_k)}$. We get
from the  obtained equation and (\ref{67}) the equation
$$
\hspace*{7ex}
\frac{1}{(x_j-x_i)^2(x_j-x_k)^2(x_j-x_l)^2}\left(-\psi_j' +
2\psi _ j\left(\frac{1}
{ x_j-x_i}+\frac{1}{x_j-x_l}+\frac{1}{x_j-x_k}\right)\right)+
$$
$$
\hspace*{7ex}
\frac{1}{(x_i-x_j)^2(x_i-x_k)^2(x_i-x_l)^2}\left(-\psi_i'+2\psi_
i\left(\frac{1}
{ x_i-x_j}+\frac{1}{x_i-x_l}+\frac{1}{x_i-x_k}\right)\right)+
$$
\begin{equation}\label{new1}
\end{equation}
$$
\hspace*{7ex}
\frac{1}{(x_k-x_i)^2(x_k-x_j)^2(x_k-x_l)^2}\left(-\psi_k'+2\psi_
k\left(\frac{1}
{ x_k-x_i}+\frac{1}{x_k-x_l}+\frac{1}{x_k-x_j}\right)\right)+
$$
$$
\hspace*{7ex}
 \frac{1}{(x_l-x_i)^2(x_l-x_k)^2(x_l-x_j)^2}\left(-\psi_l'+2\psi_
l\left(\frac{1}
{ x_l-x_i}+\frac{1}{x_l-x_j}+\frac{1}{x_l-x_k}\right)\right)=0
$$
We multiply (\ref{new1}) by $(x_l-x_j)^2$. Differentiating the obtained
equation with respect to $x_j$ we multiply the result by $(x_l-x_j)^2$. We
differentiate the obtained result with respect to $x_j$. Multiplying the
obtained
equation by $(x_l-x_j)^3(x_l-x_k)^3$, we get a polynomial of degree eight in
$x_l$. Each of its coefficients must
vanish. Equating to zero the coefficient of the highest degree, we obtain
$$
\left[\frac{1}{(x_j-x_i)^2(x_j-x_k)^2}\left(-\psi_j'+2\psi_j\left(\frac
{1}{x_j-x_i}+\frac{1}{x_j-x_k}\right)\right)\right]''+
$$
$$
\frac{1}{(x_i-x_k)^2}\left[\frac{1}{(x_i-x_j)^2}\left(-\psi_i'+2
\psi_i\left(
\frac{1}{ x_i-x_j}+\frac{1}{x_i-x_k}\right)\right)\right]''+
$$
$$
\frac{1}{(x_k-x_i)^2}\left[\frac{1}{(x_k-x_j)^2}\left(-\psi_k'+2
\psi_k\left(
\frac{1}{ x_k-x_i}+\frac{1}{x_k-x_j}\right)\right)\right]''=0
$$
We multiply the latter equation by $(x_j-x_k)^4$.
Differentiating the obtained
equation with respect to $x_j$, we multiply the result by
$(x_j-x_k)^2$ and
differentiating the obtained equation with respect to $x_j$ we obtain a
polynomial in  $x_k$.  Equating to zero the
coefficient of the highest degree, we get
$$
\left[\frac{1}{(x_j-x_i)^2}\left(-\psi_j'+\frac{2\psi_j}{x_j-x_i
}\right)\right]
^{IV}+\left[\frac{1}{(x_i-x_j)^2}\left(-\psi_i'+\frac{2\psi_i}{x
_i-x_j}\right) \right]^{IV} =0
$$
We multiply this equation by $(x_i-x_j)^6$, differentiating the obtained
equation with respect to $x_j$, we multiply the result by $(x_i-x_j)^7$.
Differentiating the obtained equation with respect to $x_j$, we get a
polynomial of $x_i$. Equating to zero the
coefficient of the highest degree we derive $\psi_j ^{VII}=0$.
Hence, the
functions $\psi_s, s=1,2,3,4$ are polynomials of degree six. Multiplying
(\ref{new1}) by $(x_j-x_i)^3$ and taking $x_j=x_i\not=x_k\not=x_l\not=x_i$,
we get from the obtained equation that $\psi_j=\psi_i$. So, we may suppose
\begin{equation}\label{d4}
\psi_s(x)=a_6x^6+a_5x^5+a_4x^4+a_3x^3+a_2x^2+a_1+a_0, \quad s=1,2,3,4.
\end{equation}
We calculate using (\ref{d4}) and (\ref{67}) that
\begin{equation}\label{d5}
R(E_i,E_j,E_j,E_i)=-\left[2(x_j-x_i)+x_k+x_l\right]\frac{a_6}{4}
-\frac{a_5}{4},
\end{equation}
$$
r_i=-a_6(\frac{3}{2}x_i+x_j+x_k+x_l)-\frac{3}{4}a_5, \quad
r_i-r_j=-\frac{a_6}{2}(x_i-x_j).
$$
Using (\ref{d5}), it is easy to verify that the conditions of
Lemma \ref{nl4} are satisfied, i.e
we have a locally conformal flat
$P$-space with four distinct Ricci eigenvalues iff $a_6\not=0$. So, we
prove VII).
\par
{\bf Spaces of type} ${\bf I_9}$. We set
$\psi_s=\frac{1}{\phi_s(x_s)}, s=2,3,4$.
Throughout the considerations in this case the indices $j,k,l$
will stand for any distinct indices
from the set $\{2,3,4\}$. We get from (\ref{d1}) and (\ref{67})
the equation $$
\frac{1}{x_j(x_j-x_k)^2(x_j-x_l)^2}\left(-\psi_j'+\psi_j\left(\frac{1}
{x_j}+2\frac{1}{x_j-x_l}+2\frac{1}{x_j-x_k}\right)\right)+
$$
\begin{equation}\label{new2}
\frac{1}{x_k(x_k-x_j)^2(x_k-x_l)^2}\left(-\psi_k'+2\psi_k\left(\frac{1}
{x_k}+2\frac{1}{x_k-x_l}+2\frac{1}{x_k-x_j}\right)\right)+
\end{equation}
$$
\frac{1}{x_l(x_l-x_k)^2(x_l-x_j)^2}\left(-\psi_l'+2\psi_l\left(\frac{1}
{x_l}+2\frac{1}{x_l-x_j}+2\frac{1}{x_l-x_k}\right)\right)=0.
$$
We multiply
(\ref{new2}) by $\frac{(x_j-x_l)^2(x_l-x_k)^2}{(x_j-x_k)}$.
Differentiating the obtained
equation with respect to $x_j$ we multiply the result by $(x_l-x_j)^2$. We
differentiate the obtained result with respect to $x_j$ and multiplying the
obtained equation by $(x_k-x_l)$, we get a polynomial of degree six in
$x_l$. Each of its coefficients must
vanish. Equating to zero the coefficient of the highest degree, we obtain
$$
\left[\frac{1}{x_j(x_j-x_k)^2}\left(-\psi_j'+\psi_j\left(\frac
{1}{x_j}+2\frac{1}{x_j-x_k}\right)\right)\right]''+
$$
$$
\left[\frac{1}{x_k(x_k-x_j)^2}\left(-\psi_k'+\psi_k\left(
\frac{1}{x_k}+2\frac{1}{x_k-x_j}\right)\right)\right]''=0
$$
Differentiating the above equation
with respect to $x_j$, we multiply the result by $(x_j-x_k)^2$. We
differentiate the obtained equation with respect to $x_j$. Multiplying the
obtained equation by $x_k$, we get a
polynomial of $x_k$ of degree five. Equating to zero the
coefficient of the highest degree we derive
$\left(\frac{\psi_j}{x_j}\right) ^{V}=0$. Hence,
the functions $\psi_s, s=2,3,4$ are polynomials of degree five having
a root equal to zero. Multiplying
(\ref{new2}) by $(x_j-x_k)^2$ and taking $x_j=x_k\not=x_l$, we get
from the obtained equation that $\psi_j=\psi_k$. So, we may suppose
\begin{equation}\label{7d4}
\psi_s(x)=a_5x^5+a_4x^4+a_3x^3+a_2x^2+a_1x, \quad s=2,3,4.
\end{equation}
We calculate using (\ref{7d4}) and (\ref{67}) that
$$
R(E_j,E_k,E_k,E_j)=-\frac{1}{4}\left(a_5\left(2(x_j+x_k)+
x_l\right)  + a_4\right),
$$
\begin{equation}\label{7d5}
R(E_1,E_j,E_j,E_1)=-\frac{1}{4}\left(a_5\left(2x_j+x_k+x_l\right)+a_4\right),
\end{equation}
$$
r_1=-a_5(x_j+x_k+x_l)-\frac{3a_4}{4},\qquad
r_j=-a_5\left(\frac{3}{2}x_j+x_k+x_l\right)-\frac{3a_4}{4}.
$$
Using (\ref{7d5}), it is easy to verify that the
conditions of Lemma \ref{nl4} are satisfied, i.e
we have a locally conformal flat
$P$-space with four distinct Ricci eigenvalues iff $a_5\not=0$.
Thus, we have proved VIII).
\par
{\bf Spaces of type} ${\bf I_{10} \quad with \quad b\not=0}$.  We set
$\psi_s=\frac{1}{\phi_s(x_s)}, s=3,4$.
Throughout the considerations in this case the indices $k,l$ will stand for
any distinct indices
from the set $\{3,4\}$. We get from (\ref{d1}) and (\ref{67}) the equation
\begin{equation}\label{new3}
\frac{1}{x_k(x_k-b)}\left(-\psi_k'+\psi_k\left(\frac{1}
{x_k}+2\frac{1}{x_k-x_l}+\frac{1}{x_k-b}\right)\right)+
\end{equation}
$$
\frac{1}{x_l(x_l-b)}\left(-\psi_l'+\psi_l\left(\frac{1}
{x_l}+2\frac{1}{x_l-x_k}+\frac{1}{x_l-b}\right)\right)=0.
$$
We differentiate (\ref{new3})  with respect to $x_k$. We
multiply the obtained result
by $(x_k-x_l)^2$. Differentiating the obtained equation by $x_k$, we get
$\left[\frac{\psi_k}{x_k(x_k-b)}\right]^{III}=0$. Multiplying
(\ref{new3}) by $(x_k-x_l)^2$ and taking $x_k=x_l$, we
derive
from the obtained equation that $\psi_3=\psi_4$. So, we may suppose
\begin{equation}\label{8d4}
\psi_s(x)=(x-b)(a_3x^3+a_2x^2+a_1x), \quad s=3,4.
\end{equation}
We calculate using (\ref{8d4}) and (\ref{67}) that
$$
R(E_1,E_2,E_2,E_1)=-\frac{1}{4}\left(a_3\left(x_3+x_4\right)+a_2\right),
$$ $$
R(E_1,E_k,E_k,E_1)=-\frac{1}{4}\left(a_3\left(2x_k+x_l\right)+a_2\right),
$$
\begin{equation}\label{8d5}
R(E_3,E_4,E_4,E_3)=-\frac{1}{4}\left(a_32\left(x_3+x_4\right)+a_
2-ba_3\right), \end{equation}
$$ \hspace*{6ex}
R(E_2,E_k,E_k,E_2)=-\frac{1}{4}\left(a_3\left(2x_k+x_l\right)+a_
2-ba_3\right), \quad
r_1=-a_3(x_k+x_k+x_l)-\frac{3a_2}{4},
$$
$$
r_2=-a_3(x_k+x_k+x_l)-\frac{3a_2}{4}+\frac{ba_3}{2},\qquad
r_k=-a_3\left(\frac{3}{2}x_k+x_j\right)-\frac{3a_2}{4}+\frac{ba_3}{2}.
$$
Using (\ref{8d5}), it is easy to verify that the conditions of
Lemma \ref{nl4}
are satisfied, i.e we have a locally conformal flat
$P$-space with four distinct Ricci eigenvalues iff $a_3\not=0$. Thua,
we get IX which completes the proof of Theorem 1.1 \hfill {\bf Q.E.D.}

\section{ Proof of Theorem \ref{1.2}}

We begin with
\begin{lem}\label{lnov}
Every locally conformal flat $4$-dimensional $\cal Q$-space is a $\cal
P$-space.
\end{lem}
{\it Proof.}
Let ({\bf M},g) be an $\cal Q$-space of dimension $4$.
Evaluating explicitly the right hand term of (\ref{1}) in  dimension $4$ we
get that the
equality (\ref{1}) is equivalent to
\begin{equation}\label{97}
   (\nabla_X
\rho)(Y,Z)=\frac{2}{9}X(s)g(Y,Z)+\frac{1}{18}Y(s)g(X,Z)+\frac{1}{18}
   Z(s)g(X,Y)
\end{equation}
for $X,Y,Z \in T_p{\bf M}, p \in {\bf M}$. Using (\ref{97}), it is easy to
verify that (\ref{10}) is satisfied. This implies that the
Jacobi operator and
its covariant derivative commute. Hence, ({\bf M},g) is a $\cal
P$-space by the results in \cite{3}. \hfill {\bf Q.E.D.}
\par
To complete the proof of Theorem 1.2 we have to pick out those
of the spaces $I)-IX)$ of Theorem
1.1 which are $\cal Q$-spaces. We keep all notations from the
previous section. Let we  set
$$  H(X,Y,Z)=\frac{2}{9}X(s)g(Y,Z)+\frac{1}{18}Y(s)g(X,Z)+\frac{1}{18}
    Z(s)g(X,Y)  $$
for $X,Y,Z \in T_p{\bf M}$, $p \in {\bf M}$. If $({\bf M},g)$ is
a space of any
of the types $I)-IX)$ then choosing $p \in \Omega$,
$U_p$
and $E_i$, $i=1,2,3,4$ as in the previous section, we obtain
that (\ref{97}) is equivalent to
\begin{equation}\label{98}
   (\nabla_i \rho)_{ii}-H_{iii}=0,
\end{equation}
\begin{equation}\label{99}
   (\nabla_i \rho)_{jj}-H_{ijj}=0,
\end{equation}
\begin{equation}\label{100}
   (\nabla_i \rho)_{ij}-H_{iij}=0,
\end{equation}
\begin{equation}\label{101}
   (\nabla_i \rho)_{jk}-H_{ijk}=0.
\end{equation}
If the Ricci tensor has not an eigenvalue of the multiplicity
three on $U_p$,
then the formulas (\ref{8}) and (\ref{9}) are equivalent to
\begin{equation}\label{102}
   E_i(r_i)=\frac{3}{2}E_i(r_k),
\end{equation}
\begin{equation}\label{103}
   E_i(r_k)=\frac{1}{4}(\nabla_k \rho)_{ki}.
\end{equation}
The equalities (\ref{102}) and (\ref{103}) together with (\ref{5}) and
(\ref{7}) imply (\ref{98}), (\ref{99}), (\ref{100})and
(\ref{101}). Thus, any locally conformal flat $\cal P$-space
which is not of
type $III_1$ is an $\cal Q$-space since  its Ricci tensor has not
three equal eigenvalues.
\par
Let ({\bf M},g) be of type $III_1)$. Then the Ricci tensor has an
eigenvalue
of multiplicity three, say $r_1=r_2=r_3\not=r_4$. The formulas
(\ref{102}) and (\ref{103}) are
valid for $i\not=4$.    The conditions (\ref{98}), (\ref{99}),
(\ref{100}) and (\ref{101}) are equivalent to
$$  (\nabla_4 \rho)_{44}-H_{444}=0,\quad (\nabla_4
\rho)_{4j}-H_{44j}=0,\quad
    (\nabla_4 \rho)_{jj}-H_{4jj}=0,\quad j \in \{1,2,3\}.  $$
Using (\ref{5}) (\ref{7}) and (\ref{8}) it is easy to see that the latter
three equalities are equivalent to the relation
$$  E_4(r_4)=\frac{3}{2}E_4(r_1).  $$
It follows from the latter equality an (\ref{19}) that the function
F=1/f is a solution of
the differential equation (\ref{w1}. This completes the proof of
Theorem 1.2. \hfill {\bf Q.E.D.}
\par
To prove the Remarks, we shall use the fact that a locally conformal flat
Riemannian $4$-manifold has parallel Ricci tensor iff it has constant Ricci
eigenvalues \cite{van,9}.\\
For Remark 2, it follows from (\ref{19}) that $\bf B^1\times_f N^3$ has
constant Ricci eigenvalues iff the function $\bf f$ on $\bf B^1$ is given by
a), b) and c). In this case (\ref{w1}) is also satisfied.
Remark 2 is contained also in \cite{Der1}.\\
For Remark 3, it follows from (\ref{24}), (\ref{cd1}) and
(\ref{32}) that  the Ricci
eigenvalues $r_1=r_2\not=r_3\not=r_4\not=r_1$ are all constant
iff the function
$f$ is a constant and $\bf B^2$ is of constant sectional curvature.\\
Remark 4 follows immediately from the formulas (\ref{nov}), (\ref{d5}),
(\ref{7d5}) and (\ref{8d5}).

\begin{flushleft}
{\bf Authors' address:}\\
         Stefan Ivanov,Irina Petrova,\\
         University of Sofia,Faculty of Math. and Inf.,\\
         bul. James Boucher 5, 1126 Sofia,\\
         BULGARIA\\
         E-mail:\quad ivanovsp@fmi.uni-sofia.bg\\
         E-mail:\quad ihp@vmei.acad.bg
\end{flushleft}


\begin{thebibliography}{9}

\bibitem {2} J.Berndt, {\em Three-dimensional Einstein-like
manifolds}, Diff. Geom. and Appl.{\bf 2}(1992), 385-397.
\bibitem {BPV} J.Berndt, F.Pr\"ufer and L.Vanhecke, {\em Symmetric-like
Riemannian manifolds and geodesic symmetries}, Proc. Royal Soc.
Edinburg A {\bf 125} (1995), 265-282.
\bibitem {btv} J.Berndt, F.Tricerri and L.Vanhecke, {\em
Generalized Heisenberg groups and Damek-Ricci harmonic soaces},
Lect. Notes Math. vol. {\bf 1598}, Springer-Verlag, 1995.
\bibitem {3} J.Berndt and L. Vanhecke, {\em Two naturally
generalizations
of locally symmetric spaces.} Diff.Geom. and Appl.{\bf 2}(1992),
57-80.
\bibitem {4} J.Berndt and L. Vanhecke, {\em Geodesic spheres
and
generalizations of symmetric spaces.} Boll.Un. Mat.Ital.A(7), {\bf 7}
(1993), no.1, 125-134.
\bibitem {5} J.Berndt and L. Vanhecke, {\em Geodesic sprays
and $\cal C$-
and $\cal P$-spaces.} Rend.Sem. Politec.Torino {\bf 50}(1992), no.4,
343-358.
\bibitem {Bes} A.Besse,  Einstein manifolds, {\em Springer,
Berlin, 1987.}
\bibitem{WB} W.Blaschke, {\em Eine Verallgemeinerung der Theorie
der confocalen $F_2$}, Math. Zeitsch. ${\bf 27} (1928), 655-668.$
\bibitem {Bou9} J.P.Bourguignon, {\em Les vari\'et\'es de
dimension 4 $\grave
a$ signature non nulle dont la courbure est harmonique sont d'Einstein},
Invent. Math. ${\bf 63} (1981), 263-286$.
\bibitem {Der1} A.Derdzinski, {\em Classification of Certain
Compact Riemannian manifolds with Harmonic Curvature and
Non-Parallel Ricci tensor}, Math. Z. ${\bf 172} (1980), 273-280$.
\bibitem {Der3} A.Derdzinski, {\em Self-dual K\"ahler manifolds and Einstein
manifolds in dimension four} Comp. Math. ${\bf 49} (1983), 405-433$.
\bibitem {De-Ch} A.Derdzinski and C.L.Chen, {\em Codazzi tensor fields,
curvature and Pontryagin form}, Proc. London Math.Soc. ${\bf 47} (1983),
15-26$.
\bibitem {e2} L.P.Eisenhart, {\em Separable systems of
St\"akel,} Ann. Math. ${\bf 35} (1934),284-305.$
\bibitem {e1} L.P.Eisenhart, {\em St\"akel systems in conformal
Euclidean space,} Ann. Math. ${\bf 36} (1935), 57-70.$
\bibitem {7} L.P.Eisenhart,  Riemannian geometry, {\em
Princeton Univ.Press, Princeton, 1949}.
\bibitem {GSV} P.Gilkey, A.Swann and L.Vanhecke, {\em
Isoparametric
geodesic spheresand a Conjecture of Osserman concerning the
Jacobi operator}, Quart.J.Math., to appear.
\bibitem {6} A.Gray, {\em Einstein-like manifolds
which are not Einstein}, Geom. Dedicata {\bf 7}(1978), 259-280.
\bibitem {8} S.Hiepko, {\em Einne innere kennzeichnung der
verzerrten produkte},  Math. Ann.{\bf 241}(1979), 209-215.
\bibitem {Hi} N.Hitchin, {\em Twistor spaces, Einstein metrics and
isomonodromic deformations}, J.Diff.Geom. ${\bf 42} N1, (1995), 30-113$.
\bibitem {9} S.Ivanov and I. Petrova. {\em Locally conformal
flat
Riemannian manifolds with constant principal Ricci curvatures and locally
conformal flat {\cal C}-spaces}, dg-ga/9702009.
\bibitem {IP1} S.Ivanov and I.Petrova, {\em Riemannian manifold
in which certain
curvature operator has constant eigenvalues along each circle},
Ann.Glob.Anal.Geom., to appear.
\bibitem {10} T.Kato, Perturbation theory for linear operators.
{\em  Springer,Berlin, 1966}.
\bibitem {1} B. O'Neil. Semi-Riemannian geometry with applications
to relativity, {\em Acad. Press, New York, 1983}.
\bibitem{Sch} J.A.Schouten.   Ricci calculus 2nd ed.,
{\em Springer 1954.}
\bibitem {11} Z.Szabo,  {\em Structure theorem on the
Riemannian space,
satisfying $R(X,Y)\circ R$, I. The local version},
J.Diff.Geom.{\bf 17}(1982), 531-582.
\bibitem{van} L.Vanhecke, Private communication.
\bibitem{Wei} Weinacht, {\em \"Uber die bedingt-periodische Bewegung eines
Massenpunktes}, Math. Ann. ${\bf 91} (1924), 279-299.$
\end{thebibliography}
\end{document}